\documentclass[10pt, a4paper]{article}

\usepackage{amsmath,amsfonts,amssymb}
\usepackage{amsthm} 
\usepackage{hyperref}
\usepackage[nameinlink,capitalize]{cleveref}
\hypersetup{colorlinks,urlcolor=blue,linkcolor = blue,citecolor = blue}
\crefname{equation}{}{} 

\usepackage{enumerate}
\usepackage{float}
\usepackage[normalem]{ulem}
\usepackage{xcolor}
\usepackage{mathrsfs}%

\usepackage{geometry}
\usepackage{bbm}
\usepackage{mathtools} 

\usepackage{enumitem}
\setlist[enumerate]{label*=\alph*),ref=\alph*)}

\newcommand{\F}{\mathscr F} 
\newcommand{\R}{\mathbb{R}}
\newcommand{\N}{\mathbb{N}}

\newtheorem{theorem}{Theorem}
\newtheorem{example}{Example}%
\newtheorem{remark}{Remark}%

\title{Correction to: Criteria for Strong and Weak Random Attractors}
\author{%
Hans~Crauel\footnote{Institute of Mathematics, Goethe-Universit\"at, Robert-Mayer-Stra{\ss}e 6-8, 60325 Frankfurt am Main, Germany. \newline E-mail: \href{mailto:crauel@math.uni-frankfurt.de}{crauel@math.uni-frankfurt.de}}
\and
Sarah~Geiss\footnote{Institute of Mathematics, TU Berlin, Stra{\ss}e des 17. Juni 136, 10623 Berlin, Germany. E-mail: \href{mailto:geiss@math.tu-berlin.de}{geiss@math.tu-berlin.de}}
\thanks{The author was supported by the DFG Research Unit FOR 2402.}
\and
Michael~Scheutzow \footnote{Institute of Mathematics, TU Berlin, Stra{\ss}e des 17. Juni 136, 10623 Berlin, Germany. E-mail:  \href{mailto:ms@math.tu-berlin.de}{ms@math.tu-berlin.de}}}

\date{September 1, 2023}   

\begin{document}

\maketitle

\begin{abstract} 
In the article 'Criteria for Strong and Weak Random Attractors' necessary and sufficient conditions for strong attractors and weak attractors are studied. In this note we correct two of its theorems on strong attractors.
\end{abstract}

[0.25em]
\textbf{MSC2020 subject classifications:}  37B25, 37C70, 37G35, 37H99, 37L55, 60D05, 60H10, 60H15, 60H25 

\bigskip


We correct two theorems which provide criteria for strong attractors given in \cite{CrauelDimitroffScheutzow}.

We use the same assumptions and notation as in \cite{CrauelDimitroffScheutzow}, i.e. let $\varphi$ be a continuous random dynamical system on a Polish space $(E,d)$ over a metric dynamical system $(\Omega, \F, (\vartheta_t)_{t\in\mathbb{R}}, P)$. We use the same letter $d$ for the complete metric on $E$ and the Hausdorff semi-distance on subsets of $E$.  For a subset $A$ of $E$ we denote the closed $\delta$-neighborhood of $A$ by $A^\delta$.

In the article the following two types of strong attractors are studied:
\begin{itemize}
\item $B$-attractors, i.e. attractors that attract all bounded subsets of $E$,
\item $C$-attractors, i.e. attractors that attract all compact subsets of $E$.
\end{itemize}

In \cite[Theorem 3.1, Theorem 3.2]{CrauelDimitroffScheutzow} the following two theorems have been stated:
\begin{theorem}[Original erroneous formulation]\label{thm:3.1}
The following are equivalent: 
\begin{enumerate}
\item $\varphi$ has a strong $B$-attractor.
\item For every $\varepsilon>0$ there exists a compact subset $C_\varepsilon$ such that for each $\delta>0$ and each bounded and closed subset $B$ of $E$ it holds that
\begin{equation*}
P\left\{ \bigcup_{s\geq 0} \bigcap_{t\geq s} \varphi(t,\vartheta_{-t} \omega) B \subseteq C_{\varepsilon}^{\delta} \right\} \geq 1 - \varepsilon.
\end{equation*}
\item There exists a compact strongly $B$-attracting set $\omega \mapsto K(\omega)$.
\end{enumerate}
\end{theorem}

\begin{theorem}[Original erroneous formulation]\label{thm:3.2}
The following are equivalent: 
\begin{enumerate}
\item $\varphi$ has a strong $C$-attractor.
\item For every $\varepsilon>0$ there exists a compact subset $C_\varepsilon$ such that for each $\delta>0$ and each compact subset $B$ of $E$ it holds that
\begin{equation*}
P\left\{ \bigcup_{s\geq 0} \bigcap_{t\geq s} \varphi(t,\vartheta_{-t} \omega) B \subseteq C_{\varepsilon}^{\delta} \right\} \geq 1 - \varepsilon.
\end{equation*}
\item There exists a compact strongly $C$-attracting set $\omega \mapsto K(\omega)$.
\end{enumerate}
\end{theorem}

The following example shows that the original formulations of Theorem \ref{thm:3.1} and Theorem \ref{thm:3.2} are incorrect.
\begin{example}
Choose $E=\R$, $\Omega=\{0\}$ and consider $\varphi(t,\omega)x := x+t$ for all $t\geq 0, x\in E$, $\omega \in \Omega$. This continuous RDS satisfies for all bounded subsets $B\subset\R$
\begin{equation*}
\begin{aligned}
\bigcup_{T\geq 0}\bigcap_{t\geq T} 
\varphi(t, \vartheta_{-t}\omega) B & \subseteq 
\bigcap_{T\geq 0}\bigcup_{t\geq T}  \varphi(t, \vartheta_{-t}\omega) B \\
& \subseteq \Omega_B(\omega) := \bigcap_{T\geq 0}\overline{\bigcup_{t\geq T} \varphi(t, \vartheta_{-t}\omega) B} = \emptyset.
\end{aligned}
\end{equation*}
This RDS has no C-attractor and hence also no B-attractor. 

In particular, (i) and (ii) of Theorem \ref{thm:3.1} and Theorem \ref{thm:3.2} of the original formulation are not equivalent. This example shows in particular that also the following stronger property is not sufficient to ensure strong $B$-attractors:
\begin{enumerate}
\item[(ii)'] For every $\varepsilon>0$ there exists a compact subset $C_\varepsilon$ such that for each $\delta>0$ and each bounded and closed subset $B$ of $E$ it holds that
\begin{equation*}
P\left\{ \Omega_B(\omega) \subseteq C^\delta_{\varepsilon} \right\} \geq 1 - \varepsilon.
\end{equation*}
\end{enumerate}
\end{example}

The following is a corrected version of Theorem \ref{thm:3.1}: The condition (ii) is modified. In addition, condition (iii) is formulated more precisely than in the original formulation.
\setcounter{theorem}{0}
\begin{theorem}[Corrected formulation]\label{thm:3.1c}
The following are equivalent: 
\begin{enumerate}
\item $\varphi$ has a strong $B$-attractor.
\item For every $\varepsilon>0$ there exists a compact subset $C_\varepsilon$ such that for each $\delta>0$ and each bounded and closed subset $B$ of $E$ there exists a $T>0$ such that 
\begin{equation*}
P\left\{ \bigcup_{t\geq T} \varphi(t,\vartheta_{-t} \omega) B \subseteq C_{\varepsilon}^{\delta} \right\} \geq 1 - \varepsilon.
\end{equation*}
\item There exists a random set $K \subseteq E\times\Omega$ such that $K(\omega)$ is $P$-a.s.\ compact and $K$ attracts all bounded subsets, i.e. 
\begin{equation*}
\lim_{t\to\infty} d(\varphi(t,\vartheta_{-t} \omega) B, K(\omega)) = 0 \quad P\text{-a.s.}
\end{equation*}
for every bounded subset $B$.
\end{enumerate}
\end{theorem}

\begin{remark}
By \cite[Lemma 3.5]{Crauel_Global} and its proof we have that 
\begin{equation*}
\bigcup_{t\geq T} \varphi(t,\vartheta_{-t} \omega) B \in \mathcal{B} \otimes \bar{\F} \quad  \text{and}  \quad  \Omega_B(\omega) \in \mathcal{B} \otimes \bar{\F}
\end{equation*}
for all bounded closed subsets $B$ of $E$. Here $\mathcal{B}$ denotes the Borel $\sigma$-algebra of $E$ and $\bar{\F}$ the $P$-completion of $\F$. Therefore, we have by the measurable projection theorem that
\begin{equation*}
\begin{aligned}
& \Omega \backslash 
\left\{ \omega \in \Omega \,\,\bigg| \,\,
\bigcup_{t\geq T} \varphi(t,\vartheta_{-t} \omega) B  \subseteq C_{\varepsilon}^\delta \right\}  \\
&  = \textrm{pr}_{\Omega}\left(
\left\{ \bigcup_{t\geq T} \varphi(t,\vartheta_{-t} \omega) B \right\} \cap \left\{ (E \backslash C^\delta_{\varepsilon})\times \Omega\right\} \right) 
 \in \bar{\F} 
 \end{aligned}
\end{equation*}
where $\textrm{pr}_{\Omega}: E \times \Omega \to \Omega$ denotes the projection onto $\Omega$. Hence, the expression in (ii) of Theorem \ref{thm:3.1c} is well-defined. 
\end{remark}

\begin{proof} 
The proof is similar to the proof presented in \cite{CrauelDimitroffScheutzow}.

Equivalence of (i) and (iii) is proven in \cite[Theorem 13]{MinimalRandomAttractors}, see also \cite[Theorem 3.4, Remark 3.5]{Crauel-7}.

We first show (i) $\implies$ (ii): Let $\varepsilon>0$ be arbitrary. Since $E$ is a Polish space and the attractor $A$ is a random variable taking values in the compact sets, there exists a compact subset $C_{\varepsilon} \subseteq E$ such that 
\begin{equation}\label{eq:1}
P\{A(\omega) \subseteq C_{\varepsilon} \} \geq 1-\varepsilon/2
\end{equation}
(see Crauel \cite[Proposition 2.15]{Crauel}). Let $B\subseteq E$ be a bounded and closed set. Then we have by (i) 
\begin{equation*}
\lim_{t\to\infty} d(\varphi(t,\vartheta_{-t} \omega) B, A(\omega)) = 0  \quad P\text{-a.s.},
\end{equation*}
i.e. for every $\delta>0$ there exists a $T(\omega)>0$ such that for all $t\geq T(\omega)$ we have 
$d(\varphi(t,\vartheta_{-t} \omega) B, A(\omega)) \leq \delta$ $P$-almost surely. Hence, there exists some deterministic $T>0$ such that 
\begin{equation}\label{eq:2}
P\left\{ \bigcup_{t\geq T} \varphi(t,\vartheta_{-t} \omega) B \subseteq A(\omega)^{\delta} \right\} \geq 1 - \varepsilon/2.
\end{equation}
Combining \eqref{eq:1} and \eqref{eq:2} implies (ii).

Now we show  (ii) $\implies$ (iii): Let $(B_k)_{k\in\N}$ be a sequence of bounded closed subsets of $E$ such that $B_0 \subseteq B_1 \subseteq B_2  \dots$ and such that for any bounded subset $B\subseteq E$ there exists some $k\in\N$ such that $B\subseteq B_k$. We modify the random attractor constructed in the proof given in \cite{CrauelDimitroffScheutzow} to ensure that it is indeed a random set: We define  $A(\omega)$ to be the  (unique) smallest closed random set that contains $\bigcup_{k\in\N} \Omega_{B_k}(\omega)$, see \cite[Proposition 17]{MinimalRandomAttractors}. By (ii) for all $\varepsilon >0$ there exists a compact set $C_{\varepsilon} \subseteq E$ such that for every $\delta >0$ and for every $k\in\N$ there exist $ T(k)>0$ such that
\begin{equation*}
P\left\{ \bigcup_{t\geq T(k)} \varphi(t,\vartheta_{-t} \omega) B_k \subseteq C_{\varepsilon}^{\delta} \right\} \geq 1 - \varepsilon.
\end{equation*}
Using that $C_{\varepsilon}$ is closed, this implies that 
$P\{\Omega_{B_k}(\omega) \subseteq C_{\varepsilon}\} \geq 1-\varepsilon$. As $\Omega_{B_k}(\omega) \subseteq \Omega_{B_{k+1}}(\omega)$ this implies  $P\{\bigcup_{k\in\N} \Omega_{B_k}(\omega) \subseteq C_{\varepsilon}\} \geq 1-\varepsilon$. This implies by the properties of $A(\omega)$ given by \cite[Proposition 17]{MinimalRandomAttractors} that $A(\omega)$ is a compact random set. 

It remains to prove that $A(\omega)$ attracts all bounded sets. To this end consider an arbitrary bounded subset $B$ of $E$ and let $k\in\N$ be such that $B\subseteq B_k$. Let $\varepsilon>0$ be arbitrary. By (ii) there exists for every $m\in\N$ some $T_m>0$ such that
\begin{equation*}
P\left\{d\left(\bigcup_{t\geq T_m} \varphi(t, \vartheta_{-t} \omega) B_k, C_{\varepsilon}\right)\leq 1/m\right\} \geq 1-\varepsilon.
\end{equation*}
which implies
\begin{equation*}
P\left\{\sup_{t\geq T_m} d(\varphi(t, \vartheta_{-t} \omega) B_k, C_{\varepsilon})\leq 1/m \text{ for infinitely many m}\right\} \geq 1-\varepsilon.
\end{equation*}
To obtain the previous inequality we used for $M_m := \{\sup_{t\geq T_m} d(\varphi(t, \vartheta_{-t} \omega) B_k, C_{\varepsilon})\leq 1/m \}$ that
\begin{equation*}
P\left[\bigcap_{n\in\N} \bigcup_{m=n}^\infty M_m \right] = \lim_{n\to\infty} P\left[\bigcup_{m=n}^\infty M_m \right] \geq \limsup_{m\to\infty} P[M_m] \geq 1 - \varepsilon.
\end{equation*}
Hence, we have
\begin{equation*}
P\left\{\lim_{t\to\infty} d(\varphi(t,\vartheta_{-t}\omega)B_k, C_\varepsilon)=0 \right\} \geq 1-\varepsilon.
\end{equation*}
Due to $\Omega_{B_k}\subseteq A$ and compactness of $C_\varepsilon$ this implies
\begin{equation*}
P[\lim_{t\to\infty} d(\varphi(t,\vartheta_{-t}\omega)B_k, A(\omega))\neq 0]<\varepsilon.
\end{equation*}
The assertion follows as the previous inequality holds for arbitrary $\varepsilon >0$.
\end{proof}
 
\newpage 

The following is a corrected version of Theorem \ref{thm:3.2}. It follows from the proof given in 
\cite{CrauelDimitroffScheutzow} and the corrected proof of Theorem \ref{thm:3.1c}. (One can use e.g. \cite[Lemma 8]{MinimalRandomAttractors} to verify that $\Omega_B$ is invariant.)
\begin{theorem}[Corrected formulation]\label{thm:3.2c}
The following are equivalent: 
\begin{enumerate}
\item $\varphi$ has a strong $C$-attractor.
\item For every $\varepsilon>0$ there exists a compact subset $C_\varepsilon$ such that for each $\delta>0$ and each compact subset $B$ of $E$ there exists a $T>0$ such that
\begin{equation*}
P\left\{ \bigcup_{t\geq T} \varphi(t,\vartheta_{-t} \omega) B \subset C_{\varepsilon}^{\delta} \right\} \geq 1 - \varepsilon.
\end{equation*}
\item There exists a random set $K \subseteq E\times\Omega$ such that $K(\omega)$ is $P$-a.s. compact and $K$ attracts all compact subsets.
\end{enumerate}
\end{theorem}


\begin{thebibliography}{1}
\bibitem{Crauel_Global}
Hans Crauel.
\newblock Global random attractors are uniquely determined by attracting
  deterministic compact sets.
\newblock {\em Ann. Mat. Pura Appl. (4)}, 176:57--72, 1999.

\bibitem{Crauel-7}
Hans Crauel.
\newblock Random point attractors versus random set attractors.
\newblock {\em J. London Math. Soc. (2)}, 63(2):413--427, 2001.

\bibitem{Crauel}
Hans Crauel.
\newblock {\em Random probability measures on {P}olish spaces}, volume~11 of
  {\em Stochastics Monographs}.
\newblock Taylor \& Francis, London, 2002.

\bibitem{CrauelDimitroffScheutzow}
Hans Crauel, Georgi Dimitroff, and Michael Scheutzow.
\newblock Criteria for strong and weak random attractors.
\newblock {\em J. Dynam. Differential Equations}, 21(2):233--247, 2009.

\bibitem{MinimalRandomAttractors}
Hans Crauel and Michael Scheutzow.
\newblock Minimal random attractors.
\newblock {\em J. Differential Equations}, 265(2):702--718, 2018.

\end{thebibliography}

\end{document}